\pdfoutput=1

\documentclass[11pt]{article}
\PassOptionsToPackage{obeyspaces}{url}
\usepackage[hidelinks]{hyperref}
\hypersetup{colorlinks = true, linkcolor = blue}
\usepackage{bookmark}
\bookmarksetup{open,numbered,depth=5}
\usepackage{amsfonts}
\usepackage{amsmath}
\usepackage{amssymb, amsthm, enumitem}
\usepackage{pgf,tikz}
\usepackage{bm}
\usepackage{cancel}
\usepackage{comment}
\usetikzlibrary{calc}

\usepackage{pgfplots}
\pgfplotsset{compat=1.13}
\usepackage{mathrsfs}
\usetikzlibrary{arrows}
\definecolor{ccqqqq}{rgb}{0.8,0,0}
\definecolor{xdxdff}{rgb}{0.49019607843137253,0.49019607843137253,1}
\definecolor{qqqqff}{rgb}{0,0,1}
\definecolor{ffqqtt}{rgb}{1,0,0.2}
\definecolor{ududff}{rgb}{0.30196078431372547,0.30196078431372547,1}
\definecolor{uuuuuu}{rgb}{0.26666666666666666,0.26666666666666666,0.26666666666666666}

\usepackage{etoolbox} 
\patchcmd{\thebibliography}{\leftmargin\labelwidth}{\leftmargin\labelwidth\addtolength\itemsep{-0.2\baselineskip}}{}{}

\usepackage{mathtools}
\usepackage{xstring}

\author{
	Ting-Wei Chao\thanks{Department of Mathematics, Massachusetts Institute of Technology, Cambridge, MA, USA. \texttt{twchao@mit.edu}} 
	\and Zichao Dong\thanks{Extremal Combinatorics and Probability Group (ECOPRO), Institute for Basic Science (IBS), Daejeon, South Korea. Supported by the Institute for Basic Science (IBS-R029-C4). Research partially supported by ERC Grant No.~882971, ``GeoScape'', in the framework of the special semester on Discrete and Convex Geometry at the Erd\H{o}s Center, Budapest, Hungary. \texttt{zichao@ibs.re.kr}. } 
	\and Zhuo Wu\thanks{Mathematics Institute and DIMAP, University of Warwick, Coventry, UK, and Extremal Combinatorics and Probability Group (ECOPRO),  Institute for Basic Science (IBS), Daejeon, South Korea. Supported by the Warwick Mathematics Institute Centre for Doctoral Training and funding from University of Warwick and the Institute for Basic Science (IBS-R029-C4). \texttt{zhuo.wu@warwick.ac.uk}. }}

\title{Empty red-red-blue triangles}
\date{}

\usepackage[nameinlink]{cleveref}

\newtheorem{theorem}{Theorem}
\newtheorem{lemma}[theorem]{Lemma}

\newtheorem{question}{Question}
\newtheorem{conjecture}{Conjecture}

\newcommand*{\eqdef}{\stackrel{\mbox{\normalfont\tiny def}}{=}} 
\newcommand*{\veps}{\varepsilon}                                
\newcommand{\Yemph}[1]{\textcolor{magenta}{\emph{#1}}}
\definecolor{darkpastelgreen}{rgb}{0.01, 0.75, 0.24}

\newcommand*{\R}{\mathbb{R}}                                    
\newcommand*{\C}{\mathbb{C}}                                    
\newcommand*{\cC}{\mathcal{C}}
\newcommand*{\cS}{\mathcal{S}}
\newcommand*{\cH}{\mathcal{H}}
\newcommand*{\cG}{\mathcal{G}}
\newcommand*{\cT}{\mathcal{T}}
\newcommand*{\Nrrb}{\mathcal{N}_{\mathsf{r}\mathsf{r}\mathsf{b}}}

\DeclareMathOperator{\conv}{conv}
\DeclareMathOperator{\disc}{disc}

\newcommand{\lc}{\left\lceil}
\newcommand{\rc}{\right\rceil}

\allowdisplaybreaks[4]

\begin{document}
	
	\maketitle
	
	\begin{abstract}
		Let $P$ be a $2n$-point set in the plane that is in general position. We prove that every red-blue bipartition of $P$ into $R$ and $B$ with $|R| = |B| = n$ generates $\Omega(n^{3/2})$ red-red-blue empty triangles. 
	\end{abstract}
	
	\section{Introduction}
	
	\paragraph{The study of convex holes.} A set of points in the plane is \Yemph{in general position} if no three points from the set are collinear, and a set of points in the plane is \Yemph{in convex position} if it serves as the vertices of a convex polygon. A classical result due to Erd\H{o}s and Szekeres \cite{erdos_szekeres} states that every set of $4^k$ points in the plane that is in general position contains a $k$-subset that is in convex position. 
	
	Let $X \subset \R^2$ be a finite point set that is in general position. Then $Y \subseteq X$ is a \Yemph{$k$-hole} if $|Y| = k$, $Y$ is in convex position and the convex hull of $Y$ contains no point of $X$ other than $Y$. Erd\H{o}s raised the question on whether a sufficiently large $X$ (i.e.,~$|X|$ is large with respect to $k$) always contains a $k$-hole. This is easily seen to be true for $k = 3, 4, 5$. Astonishingly, Horton \cite{horton} disproved the existence of $7$-holes by explicit constructions (later known as \emph{Horton sets}). Many years later, Nicol\'{a}s \cite{nicolas} and Gerken \cite{gerken} independently established the existence of $6$-holes. 
	
	\paragraph{Bicolored empty triangles.} The existence of a $3$-hole in an $n$-point set is trivial (as seen by the triangle of smallest area), and a further problem is to determine the smallest number of $3$-holes in an $n$-set. The best-known lower bound is $n^2 + \Omega(n \log^{2/3} n)$ (due to Aichholzer et al.~\cite{aichholzer_etal_5hole}), while the current upper bound is roughly $1.6196n^2$ (due to B\'{a}r\'{a}ny and Valtr \cite{barany_valtr}). 
	
	In this note, we consider the minimum number of $3$-holes under a bicolored setting. Let $R \subset \R^2$ be a set of finitely many red points, and $B \subset \R^2$ be a set of finitely many blue points. Moreover, we assume that $R \cup B$ is in general position. 
	
	We say that $a, b, c \in R \cup B$ form an \Yemph{empty triangle}, if $\conv(\{a, b, c\}) \cap (R \cup B) = \{a, b, c\}$. Here we denote by $\conv(P)$ the convex hull of any finite point set $P$. Moreover, 
	\begin{itemize}
		\vspace{-0.5em}
		\item if $|\{a, b, c\} \cap R| = 3$ and $|\{a, b, c\} \cap B| = 0$, then they form an empty \emph{red}-\emph{red}-\emph{red} triangle; 
		\vspace{-0.5em}
		\item if $|\{a, b, c\} \cap R| = 2$ and $|\{a, b, c\} \cap B| = 1$, then they form an empty \emph{red}-\emph{red}-\emph{blue} triangle; 
		\vspace{-0.5em}
		\item if $|\{a, b, c\} \cap R| = 1$ and $|\{a, b, c\} \cap B| = 2$, then they form an empty \emph{red}-\emph{blue}-\emph{blue} triangle; 
		\vspace{-0.5em}
		\item if $|\{a, b, c\} \cap R| = 0$ and $|\{a, b, c\} \cap B| = 3$, then they form an empty \emph{blue}-\emph{blue}-\emph{blue} triangle. 
	\end{itemize}
	\vspace{-0.25em}
	
	Suppose $R$ and $B$ are both $n$-point sets in the plane. We make the following direct observations concerning the number of empty triangles in red and blue: 
	\vspace{-0.5em}
	\begin{itemize}
		\item An empty red-red-red does not necessarily exist. To see this, we begin by choosing $\veps > 0$ to be sufficiently small in terms of $n$. Identify the real plane $\R^2$ and the complex plane $\C$. Put $n$ red points at $\exp(\frac{2\pi\text{i}}{n})$ and $n$ blue points at $(1 - \veps)\exp(\frac{2\pi\text{i}}{n})$, respectively. Then it is not difficult to see that no empty red-red-red triangle exists in such a configuration. 
		\vspace{-0.5em}
		\item The number of empty red-red-blue triangles and empty red-blue-blue triangles adds up to some number that is quadratic in $n$. In fact, each pair of red and blue points $p, q$ is contained in at least one empty red-red-blue or empty red-blue-blue triangle, since the point (other than $p, q$) that is closest to the line $\overline{pq}$ forms such an empty triangle with $p, q$. 
	\end{itemize}
	\vspace{-0.5em}
	These observations lead to the following pair of questions: 
	
	\begin{question} \label{ques:mono}
		What is the minimum number of empty monochromatic triangles in $R \cup B$? 
	\end{question}
	
	\begin{question} \label{ques:RRB}
		What is the minimum number of empty red-red-blue triangles in $R \cup B$? 
	\end{question}
	
	We remark that monochromatic refers to red-red-red or blue-blue-blue. Each  of \Cref{ques:mono} and \Cref{ques:RRB} has an obvious $O(n^2)$ upper bound (since there exists a set of $n$ points containing $O(n^2)$ many empty triangles \cite{barany_valtr}). These quadratic upper bounds seems close to the truth, but proving a superlinear lower bound to any of the questions is already non-trivial. 
	
	Aichholzer et al.~\cite{aichholzer_etal} initiated the study of \Cref{ques:mono} and proved a $\Omega(n^{5/4})$ lower bound. This bound was later improved to $\Omega(n^{4/3})$ by Pach and T\'{o}th \cite{pach_toth}. 
	
	In this note, we make some progress on \Cref{ques:RRB}. To be specific, we establish a $\Omega(n^{3/2})$ lower bound on the number of red-red-blue triangles in a strong sense (see \Cref{thm:RRB_3/2}). 
	
	As a curious reader may wonder, what if we ask similar questions concerning more color classes? 
	\vspace{-0.5em}
	\begin{itemize}
		\item Devillers et al.~\cite{devillers_etal} proved that for any integer $n$, there exist $n$-sets $R, G, B$ (with $R \cup G \cup B$ being in general position) containing no monochromatic empty triangle.
		\vspace{-0.5em}
		\item Fabila-Monroy et al.~\cite{fabila-monroy_etal} showed that for every integer $k \ge 3$, the minimum number of empty rainbow triangles (i.e., triangles whose vertices are of pairwise different colors) in an $k$-colored $n$-point set $X$ is of order $\Theta \bigl( k^2 \cdot \min (k, n) \bigr)$, which is independent on $n$ as $n \to \infty$. 
	\end{itemize}
	\vspace{-0.5em}
	These results suggest that our questions concerning two colors are particularly interesting. Other related problems concerning colored point sets can be found in the survey papers \cite{kaneko_kano,kano_urrutia}. For a extensive list of related combinatorial geometry problems, we refer the readers to \cite{brass_moser_pach}. 
	
	\paragraph{Our result.} Here and after we fix integers $m \ge n \ge 2$ and assume that $|R| = n, \, |B| = m$. This setting is more general than the previous ``$|R| = |B| = n$'' one, and so our \Cref{thm:RRB_3/2} below implies a $\Omega(n^{3/2})$ lower bound towards \Cref{ques:RRB}. 
	
	\begin{theorem} \label{thm:RRB_3/2}
		There exist $\Omega(n^{3/2})$ empty red-red-blue triangles in $R \cup B$. 
	\end{theorem}
	
	\section{Proof of \texorpdfstring{\Cref{thm:RRB_3/2}}{Theorem 1}}
	
	To prove a lower bound on $\Nrrb$, the idea is to analyze the discrepancy between red and blue points of a region on the plane. Similar ideas were already applied before. See for instance \cite[Lemma 3]{aichholzer_etal_count}. 
	
	Let $\cC$ be a convex region (i.e.~a convex subset of $\R^2$). Define its \Yemph{discrepancy} as the quantity
	\[
	\disc(\cC) \eqdef |R(\cC)| - |B(\cC)|, \qquad \text{where $R(\cC) \eqdef \cC \cap R$ and $B(\cC) \eqdef \cC \cap B$}. 
	\]
	Keep in mind that the discrepancy of $\cC$ is negative if it contains more blue points than red points. Denote by $\Nrrb(\cC)$ the number of empty red-red-blue triangles formed by points in $\cC$. 
	
	\newpage
	
	\begin{lemma} \label{lem:discrepancy}
		We have $\Nrrb(\cC) \ge |B(\cC)| \cdot \disc(\cC)$ for any convex region $\cC$. 
	\end{lemma}
	
	\begin{proof}
		Assume that $|B(\cC)| > 0$ and $\disc(\cC) > 0$, for otherwise the statement is trivially true. Fix an arbitrary blue point $b \in B(\cC)$. We claim that $b$ contributes at least $\disc(\cC)$ empty red-red-blue triangles in $\cC$. Then the lemma is proved by summing over $b$. 
		
		To see this, we draw a ray from $b$ to every red point in $R(\cC)$. These $|R(\cC)|$ rays partition $\cC$ into $|R(\cC)|$ sectors. In each sector containing no blue point (later called \Yemph{empty sector}), $b$ forms an empty red-red-blue triangle with the two red points on the boundary of the sector. However, there is at most one exception, which is the unique (if exists) sector whose angle at the apex is bigger than $\pi$. Since there are exactly $|B(\cC)| - 1$ blue points other than $b$ in $\cC$, the number of empty sectors is bounded from below by $|R(\cC)| - \bigl( |B(\cC)| - 1 \bigr) = \disc(\cC) + 1$, and so the claim is true. 
	\end{proof}
	
	For any red point $r \in R$, we are going to partition $\R^2$ into sectors sharing the apex $r$. Here we abuse the terminology ``partition'' because indeed these sectors can have boundary rays (including $r$) in common. Draw a ray from $r$ to every other red point. Such rays partition $\R^2$ into $n - 1$ sectors. For each sector containing at least one blue point, we call it a \Yemph{blue sector with respect to $r$}. Denote by $p(r)$ the number of blue sectors with respect to $r$. For each ray that is shared by two adjacent non-blue sectors, we erase it and hence merge the adjacent sectors. This erasing leaves us exactly $p(r)$ \Yemph{red sectors with respect to $r$}. See \hyperlink{figone}{Figure~1} for an illustration. Define $p \eqdef \min\limits_{r \in R} p(r)$. 
	
	\vspace{0.75em}
	\begin{center}
		\begin{tikzpicture}[line cap=round,line join=round,>=triangle 45,x=1.5cm,y=1.5cm, scale = 0.9]
			\clip (-2.5,-2) rectangle (2.5,2);
			\fill[line width=0pt,color=qqqqff,fill=qqqqff,fill opacity=0.1] (0,0) -- (3.158743596232484,4.352104629662982) -- (-1.8688361505976068,4.309877964871542) -- cycle;
			\fill[line width=0pt,color=ccqqqq,fill=ccqqqq,fill opacity=0.1] (0,0) -- (-1.8688361505976068,4.309877964871542) -- (-3.3807033920851097,3.4383408838173914) -- cycle;
			\fill[line width=0pt,color=qqqqff,fill=qqqqff,fill opacity=0.1] (0,0) -- (-3.3807033920851097,3.4383408838173914) -- (-4.286904663758825,1.2901071025651283) -- cycle;
			\fill[line width=0pt,color=qqqqff,fill=qqqqff,fill opacity=0.1] (0,0) -- (-4.286904663758825,1.2901071025651283) -- (-3.6,-0.4) -- cycle;
			\fill[line width=0pt,color=ccqqqq,fill=ccqqqq,fill opacity=0.1] (0,0) -- (-3.6,-0.4) -- (-4.015540973074522,-1.3186043720619924) -- cycle;
			
			\fill[line width=0pt,color=qqqqff,fill=qqqqff,fill opacity=0.1] (0,0) -- (-4.015540973074522,-1.3186043720619924) -- (-1.8387502236523767,-3.241327126766) -- cycle;
			\fill[line width=0pt,color=ccqqqq,fill=ccqqqq,fill opacity=0.1] (0,0) -- (-1.8387502236523767,-3.241327126766) -- (2.5159023479201963,-3.1840649061015043) -- cycle;
			\fill[line width=0pt,color=qqqqff,fill=qqqqff,fill opacity=0.1] (0,0) -- (2.5159023479201963,-3.1840649061015043) -- (4.006417767326709,-1.4630083744286697) -- cycle;
			\fill[line width=0pt,color=ccqqqq,fill=ccqqqq,fill opacity=0.1] (0,0) -- (4.006417767326709,-1.4630083744286697) -- (4.2031718234482955,0) -- cycle;
			\fill[line width=0pt,color=qqqqff,fill=qqqqff,fill opacity=0.1] (0,0) -- (4.1144848763510495,2.3060614580075653) -- (4.2031718234482955,0) -- cycle;
			\fill[line width=0pt,color=ccqqqq,fill=ccqqqq,fill opacity=0.1] (0,0) -- (4.1144848763510495,2.3060614580075653) -- (3.158743596232484,4.352104629662982) -- cycle;
			\draw [line width=0.8pt,domain=0:5.322634841763504] plot(\x,{(-0--1.2968034849288912*\x)/0.9412158144525837});
			\draw [line width=0.8pt,domain=-3.4223320580554053:0] plot(\x,{(-0--1.5050905149313933*\x)/-0.6526327629578696});
			\draw [line width=0.8pt,domain=-3.4223320580554053:0] plot(\x,{(-0--1.1519081597097593*\x)/-1.1325985277462447});
			\draw [line width=0.8pt,domain=-3.4223320580554053:0] plot(\x,{(-0--0.5179911118760573*\x)/-1.7212357864489691});
			\draw [line width=0.8pt,domain=-3.4223320580554053:0] plot(\x,{(-0-0.532499995962649*\x)/-1.6216202503608157});
			\draw [line width=0.8pt,domain=-3.4223320580554053:0] plot(\x,{(-0-1.166417043796351*\x)/-0.6616887207840654});
			\draw [line width=0.8pt,domain=0:5.322634841763504] plot(\x,{(-0-1.202640875101134*\x)/0.9502717722787795});
			\draw [line width=0.8pt,domain=0:5.322634841763504] plot(\x,{(-0-0.6049476585722149*\x)/1.6566364827220486});
			\draw [line width=0.8pt,domain=0:5.322634841763504] plot(\x,{(-0-0*\x)/2});
			\draw [line width=0.8pt,domain=0:5.322634841763504] plot(\x,{(-0--0.8168377201405169*\x)/1.457405410545742});
			\draw [line width=0.8pt,domain=0:-5.322634841763504] plot(\x,{(-0-0.4*\x)/-3.6});
			\draw [fill=ffqqtt, draw=none] (0,0) circle (2pt);
			\draw [fill=ffqqtt, draw=none] (2,0) circle (2pt);
			\draw [fill=ffqqtt, draw=none] (1.457405410545742,0.8168377201405169) circle (2pt);
			\draw [fill=ffqqtt, draw=none] (0.9412158144525837,1.2968034849288912) circle (2pt);
			\draw [fill=ffqqtt, draw=none] (-0.6526327629578696,1.5050905149313933) circle (2pt);
			\draw [fill=ffqqtt, draw=none] (-1.1325985277462447,1.1519081597097593) circle (2pt);
			\draw [fill=ffqqtt, draw=none] (-1.7212357864489691,0.5179911118760573) circle (2pt);
			\draw [fill=ffqqtt, draw=none] (-1.6216202503608157,-0.532499995962649) circle (2pt);
			\draw [fill=ffqqtt, draw=none] (-0.6616887207840654,-1.166417043796351) circle (2pt);
			\draw [fill=ffqqtt, draw=none] (0.9502717722787795,-1.202640875101134) circle (2pt);
			\draw [fill=ffqqtt, draw=none] (1.6566364827220486,-0.6049476585722149) circle (2pt);
			\draw [fill=qqqqff, draw=none] (1.0589432661931286,0.27348025056877223) circle (2pt);
			\draw [fill=qqqqff, draw=none] (1.339677958805197,0.49987919622366583) circle (2pt);
			\draw [fill=qqqqff, draw=none] (1.5570209466338951,0.21008854578540204) circle (2pt);
			\draw [fill=qqqqff, draw=none] (0.33446664009746796,0.9073972984024743) circle (2pt);
			\draw [fill=qqqqff, draw=none] (0.2801308931402934,1.3149154005812829) circle (2pt);
			\draw [fill=qqqqff, draw=none] (-0.24511466077906055,1.1519081597097593) circle (2pt);
			\draw [fill=ffqqtt, draw=none] (-0.7341363833936314,1.0341807079692147) circle (2pt);
			\draw [fill=qqqqff, draw=none] (-0.946158634809191,0.5527636685160655) circle (2pt);
			\draw [fill=qqqqff, draw=none] (-0.6092109854434086,0.3803718479103166) circle (2pt);
			\draw [fill=qqqqff, draw=none] (-0.824695961655589,-0.6230595742246063) circle (2pt);
			\draw [fill=qqqqff, draw=none] (-0.48056956426015024,-0.36043679726492983) circle (2pt);
			\draw [fill=ffqqtt, draw=none] (0.1895713148783358,-0.5506119116150404) circle (2pt);
			\draw [fill=ffqqtt, draw=none] (-0.2179467873004733,-0.885682351184283) circle (2pt);
			\draw [fill=qqqqff, draw=none] (0.9593277301049752,-0.7045631946603681) circle (2pt);
			\draw [fill=ffqqtt, draw=none] (0.8506562361906261,0.8983413405762786) circle (2pt);
			\draw [fill=ffqqtt, draw=none] (1.113279013150303,0.8077817623143211) circle (2pt);
			\draw [fill=qqqqff, draw=none] (2.1,0.5) circle (2pt);
			\draw [fill=qqqqff, draw=none] (2.2,0.8) circle (2pt);
			\draw [fill=qqqqff, draw=none] (0,0.5) circle (2pt);
			\draw [fill=qqqqff, draw=none] (0.1,1.8) circle (2pt);
			\draw [fill=qqqqff, draw=none] (0.4,1.7) circle (2pt);
			\draw [fill=qqqqff, draw=none] (-2.1,1.25) circle (2pt);
			\draw [fill=qqqqff, draw=none] (-1.6,1) circle (2pt);
			\draw [fill=qqqqff, draw=none] (-1.55,-1.05) circle (2pt);
			\draw [fill=qqqqff, draw=none] (-1.25,-1.2) circle (2pt);
			\draw [fill=qqqqff, draw=none] (2.1,-0.98) circle (2pt);
			\draw [fill=qqqqff, draw=none] (1.55,-1.6) circle (2pt);
			\draw [fill=qqqqff, draw=none] (-2.2,0.3) circle (2pt);
			\draw [fill=ffqqtt, draw=none] (-1.8,-0.2) circle (2pt);
			\node at (0.02, -0.24) {$r$};
		\end{tikzpicture}
		
		\begin{tikzpicture}
			\clip (-5, -0.25) rectangle (5, 0.75);
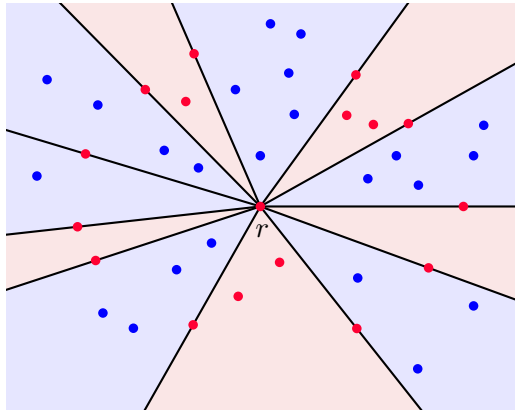
			\node at (0, 0) {\textbf{\hypertarget{figone}{Figure 1}:} Blue sectors and red sectors with respect to $r$. };
		\end{tikzpicture}
	\end{center}
	
	We remark that any sector refers to a closed region sandwiched between two rays. The two rays bounding a red sector possibly coincide when the sector contains only one red point other than $r$. 
	
	\newpage
	
	\begin{lemma} \label{lem:red_center}
		$\Nrrb(\R^2) \ge pn/2$. 
	\end{lemma}
	
	\begin{proof}
		Fix an arbitrary red point $r \in R$ and an arbitrary blue sector $\cS$ with respect to $r$. Assume $\cS$ is the sector between rays $\overrightarrow{rr_1}$ and $\overrightarrow{rr_2}$, where $r_1, r_2 \in R$. Pick $b \in \cS \cap B$ so that the smaller one of the angles $\angle \bigl( \overrightarrow{\mathstrut rb}, \overrightarrow{\mathstrut rr_1} \bigr), \angle \bigl( \overrightarrow{\mathstrut rb}, \overrightarrow{\mathstrut rr_2} \bigr)$ is minimized. Note that one of the angles is possibly bigger than $\pi$, but the smaller one cannot. Assume without loss of generality that the minimum is achieved at $\angle \bigl( \overrightarrow{\mathstrut rb}, \overrightarrow{\mathstrut rr_1} \bigr)$. It follows that $r, r_1, b$ form an empty red-red-blue triangle. By summing over all blue sectors with respect to $r$ and all red points $r$, we obtain at least $pn$ empty red-red-blue triangles. Since each such triangle is counted at most twice (from its two red vertices), we obtain $\Nrrb(\R^2) \ge pn/2$. 
	\end{proof}
	
	\begin{proof}[Proof of \Cref{thm:RRB_3/2}]
		\Cref{lem:red_center} implies that we are done if $p > \sqrt{n}$. Assume $p \le \sqrt{n}$ then, and fix a red point $r_0 \in R$ with $p(r_0) = p$. For convenience purpose, we suppose that $n \ge 5$. 
		
		Draw a line $\ell$ through $r_0$ that bisects $R$. Formally, we fix a line $\ell$ through $r_0$ that separates the plane into two closed half-planes $\cH_+, \cH_-$ with $|R(\cH_+)| = |R(\cH_-)| = \bigl\lceil (n+1)/2 \bigr\rceil \ge 3$. Depending on the parity of $n$, our $\ell$ goes through $0$ or $1$ point of $R$ other than $r_0$. Due to the discrete intermediate value theorem, the existence of $\ell$ is obvious. By applying a possible perturbation, we assume further that $\ell$ goes through no point of $B$. Introduce coordinates so that $r_0$ is the origin and $\ell$ is the $x$-axis. Choose the direction of the $x$-axis appropriately so that $\cH \eqdef \bigl\{ (x, y) \in \R^2 : y \ge 0 \bigr\}$ contains at least a half of the blue points (i.e., $|B(\cH)| \ge \lc m/2 \rc$). Observe that $|B(\cH)| \ge |R(\cH)| - 1$. 
		
		Consider the following algorithm: At the beginning, a ray $\bm{v}$ starts from the positive half of the $x$-axis $\bigl\{ (x, 0) : x \ge 0 \bigr\}$ and rotates counterclockwise around the origin $r_0$ to scan through $\cH$. This ray will generate a collection of good sectors, where a \Yemph{good sector} is a sector $\cG$ with $|B(\cG)| \ge |R(\cG)|/3$. 
		\vspace{-0.25em}
		\begin{itemize}[leftmargin = 1.8cm]
			\item[\textsc{Step $i$}:] After \textsc{Step $i-1$} (the beginning if $i = 1$), when the ray $\bm{v}$ hits the first red point, we label $\bm{v}$ at this moment as $\bm{v}_i^-$. At the first time when $\bm{v}$ hits some blue point so that the number of blue points is bigger than or equal to a third of the number of red points in the sector between $\bm{v}_i^-$ and $\bm{v}$, we label $\bm{v}$ at this moment as $\bm{v}_i^+$. If both $\bm{v}_i^-$ and $\bm{v}_i^+$ exist, then $\cG_i$ is defined as the sector between them. If $\bm{v}$ arrives the negative half of the $x$-axis before one of $\bm{v}_i^-, \bm{v}_i^+$ appears, then $\cG_i$ is not generated and the algorithm terminates. 
		\end{itemize}
		\vspace{-0.25em}
		Let $\cG_1, \dots, \cG_k$ (each containing $r_0$) be all good sectors generated by the algorithm. This means that the algorithm terminates in \textsc{Step $k+1$}. Denote by $\cT$ the sector between $\bm{v}_{k+1}^-$ and the negative half of the $x$-axis. If no red points other than $r_0$ exist after $\cG_k$, then $\cT$ denotes the negative half of the $x$-axis (hence $|R(\cT)| = 1$ and $|B(\cT)| = 0$). See \hyperlink{figtwo}{Figure~2} for an illustration. 
		
		\vspace{0.5em}
		\begin{center}
			\begin{tikzpicture}[line cap=round,line join=round,>=triangle 45,x=2cm,y=2cm, scale = 0.9]
				\clip(-2,-0.7) rectangle (2,2);
				\fill[line width=0pt,color=qqqqff,fill=qqqqff,fill opacity=0.1] (0,0) -- (4.617395575649076,0) -- (4.492019706117622,1.2227409713247717) -- cycle;
				\fill[line width=0pt,color=ccqqqq,fill=ccqqqq,fill opacity=0.1] (0,0) -- (4.492019706117622,1.2227409713247717) -- (4.405823795814748,2.398139748182151) -- cycle;
				\fill[line width=0pt,color=ccqqqq,fill=ccqqqq,fill opacity=0.1] (0,0) -- (2.1725661197857247,2.954495169227977) -- (0.02166145868131356,3.2899724801546206) -- cycle;
				\fill[line width=0pt,color=qqqqff,fill=qqqqff,fill opacity=0.1] (0,0) -- (-0.8611262409352076,3.2290905698362398) -- (-1.3938429562210393,3.046444838881098) -- cycle;
				\fill[line width=0pt,color=qqqqff,fill=qqqqff,fill opacity=0.1] (0,0) -- (-3.735771731884041,1.324608865319078) -- (-1.3938429562210393,3.046444838881098) -- cycle;
				\fill[line width=0pt,color=qqqqff,fill=qqqqff,fill opacity=0.1] (0,0) -- (-0.8611262409352076,3.2290905698362398) -- (0.02166145868131356,3.2899724801546206) -- cycle;
				\fill[line width=0pt,color=qqqqff,fill=qqqqff,fill opacity=0.1] (0,0) -- (2.1725661197857247,2.954495169227977) -- (4.405823795814748,2.398139748182151) -- cycle;
				\fill[line width=0pt,color=ccqqqq,fill=ccqqqq,fill opacity=0.1] (0,0) -- (-3.4771840009754174,0) -- (-3.735771731884041,1.324608865319078) -- cycle;
				
				\draw [line width=0.8pt,domain=-4.026985577491008:4.46604091192311,dashed] plot(\x,{(-0-0*\x)/0.5});
				\draw [line width=0.8pt,domain=0:4.46604091192311] plot(\x,{(-0--1.2227409713247717*\x)/4.492019706117622});
				\draw [line width=0.8pt,domain=0:4.46604091192311] plot(\x,{(-0--2.954495169227977*\x)/2.1725661197857247});
				\draw [line width=0.8pt,domain=-4.026985577491008:0] plot(\x,{(-0--1.324608865319078*\x)/-3.735771731884041});
				\draw [line width=0.8pt,domain=-4.026985577491008:0,color=ccqqqq,opacity=0.1] plot(\x,{(-0--3.046444838881098*\x)/-1.3938429562210393});
				\draw [line width=0.8pt,domain=0:4.439059156213827] plot(\x,{(-0--3.3013878383393167*\x)/3.598473689886184});
				\draw [line width=0.8pt,domain=-4.026985577491008:0] plot(\x,{(-0--3.251921286205633*\x)/-2.718024505645821});
				
				\draw [line width=0.8pt,domain=0:0.2] plot({0.8066551545160091+\x*(cos(0))},{0.7400586321744388+\x*(sin(0))});
				\draw [line width=0.8pt,domain=0:0.2] plot({0.6485206090503868+\x*(cos(125))},{0.8819298934723003+\x*(sin(125))});
				\draw [line width=0.8pt,domain=0:0.2] plot({-1.0317664160320579+\x*(cos(200))},{0.36583791508195346+\x*(sin(200))});
				\draw [line width=0.8pt,domain=0:0.2] plot({-0.7020456059686381+\x*(cos(60))},{0.8399471914967488+\x*(sin(60))});
				\draw [line width=0.8pt,domain=0:0.2] plot({-1.0947051280427746+\x*(cos(135))},{0+\x*(sin(135))});
				\draw [line width=0.8pt,domain=0:0.2] plot({1.0562723662958748+\x*(cos(50))},{0.2875204437302872+\x*(sin(50))});
				
				\draw [shift={(0,0)},line width=0.8pt]  plot[domain=0.9367421479549861:2.2670010123003084,variable=\t]({1*1.0947051280427746*cos(\t r)+0*1.0947051280427746*sin(\t r)},{0*1.0947051280427746*cos(\t r)+1*1.0947051280427746*sin(\t r)});
				\draw [shift={(0,0)},line width=0.8pt]  plot[domain=2.8008485084552643:3.141592653589793,variable=\t]({1*1.0947051280427749*cos(\t r)+0*1.0947051280427749*sin(\t r)},{0*1.0947051280427749*cos(\t r)+1*1.0947051280427749*sin(\t r)});
				\draw [shift={(0,0)},line width=0.8pt]  plot[domain=0.265763951616721:0.7423679567133503,variable=\t]({1*1.0947051280427749*cos(\t r)+0*1.0947051280427749*sin(\t r)},{0*1.0947051280427749*cos(\t r)+1*1.0947051280427749*sin(\t r)});
				
				\draw [fill=ffqqtt, draw=none] (0,0) circle (2pt);
				
				\draw[color=black] (0,-0.15) node {\footnotesize $r_0$};
				\draw[color=black] (-0.05,1.22) node {\footnotesize $\cG_2$};
				\draw[color=black] (-1.32,0.22) node {\footnotesize $\cT$};
				\draw[color=black] (1.1,0.6) node {\footnotesize $\cG_1$};
				
				\draw [fill=ffqqtt, draw=none] (-1.5917091647557768,0) circle (2pt);
				\draw [fill=ffqqtt, draw=none] (-0.778822214656253,0.2761503871439254) circle (2pt);
				\draw [fill=ffqqtt, draw=none] (-0.7501682436628383,1.6396009062564258) circle (2pt);
				\draw [fill=ffqqtt, draw=none] (0.004752251305354271,0.7217785396364786) circle (2pt);
				\draw [fill=ffqqtt, draw=none] (1.0619597233890699,1.4441700273670546) circle (2pt);
				\draw [fill=ffqqtt, draw=none] (1.599970966599487,0.8708823023255056) circle (2pt);
				\draw [fill=ffqqtt, draw=none] (1.7676837407854369,0.48116871153957375) circle (2pt);
				\draw [fill=ffqqtt, draw=none] (-0.8687364797250052,0.12030802420392708) circle (2pt);
				\draw [fill=ffqqtt, draw=none] (-1.8200163284497048,0.32578447152846185) circle (2pt);
				\draw [fill=qqqqff, draw=none] (-1.4740963834589054,1.2950521573926819) circle (2pt);
				\draw [fill=qqqqff, draw=none] (-1.667811552653753,0.8661114256040903) circle (2pt);
				\draw [fill=qqqqff, draw=none] (-1.3938429562210393,1.4216588572593138) circle (2pt);
				\draw [fill=qqqqff, draw=none] (-0.6632600324004702,0.6758554558591506) circle (2pt);
				\draw [fill=qqqqff, draw=none] (-0.9905003003617667,0.9269933359224709) circle (2pt);
				\draw [fill=ffqqtt, draw=none] (0.32953020972312486,1.0058630833803737) circle (2pt);
				\draw [fill=ffqqtt, draw=none] (0.6034988061558383,1.622292425353978) circle (2pt);
				\draw [fill=ffqqtt, draw=none] (0.7902955764508708,0.341004949108057) circle (2pt);
				\draw [fill=qqqqff, draw=none] (1.5547786548805378,0.1915675328720319) circle (2pt);
				\draw [fill=qqqqff, draw=none] (-0.22921946580604066,1.878304542298994) circle (2pt);
				\draw [fill=qqqqff, draw=none] (1.8171696698845659,1.2841018685435752) circle (2pt);
				\draw [fill=qqqqff, draw=none] (1.5020203663260405,1.3780152913897212) circle (2pt);
				\draw [fill=qqqqff, draw=none] (-0.9770850399750475,1.169012138532197) circle (2pt);
				\draw [fill=qqqqff, draw=none] (-1.3758796755478557,1.1165539397758957) circle (2pt);
				\draw [fill=qqqqff, draw=none] (-0.5,0.5) circle (2pt);
				\draw [fill=qqqqff, draw=none] (-1.1772970706567991,1.2798316798130869) circle (2pt);
				\draw [fill=qqqqff, draw=none] (-1.6186909204650597,1.4624774107682288) circle (2pt);
				\draw [fill=qqqqff, draw=none] (-1.519757816197691,1.7212255296213468) circle (2pt);
				\draw [fill=qqqqff, draw=none] (-1.801336651420202,1.6679538580927638) circle (2pt);
				\draw [fill=ffqqtt, draw=none] (1.3873534015049906,0.5644692335721136) circle (2pt);
				
				\node at (0, -0.55) {\textbf{\hypertarget{figtwo}{Figure 2}:} Generating $\cG_1, \dots, \cG_k$ and $\cT$. };
				
			\end{tikzpicture}
		\end{center}
		
		To finish the proof, we need to analyze some detailed properties of our algorithm: 
		\vspace{-0.5em}
		\begin{itemize}
			\item Each $\bm{v}_i^-$ comes from a distinct red sector with respect to $r_0$, and so
			\begin{equation} \label{eq:gs_stepnumber}
				k \le p \le \sqrt{n}. 
			\end{equation}
			\item The way we defined $\cG_i$ implies that $\bm{v}_i^- \in \cG_i, \, r_0 \in \cG_i$, and $|B(\cG_i)| = \bigl\lceil |R(\cG_i)|/3 \bigr\rceil$. Then
			\[
			\disc(\cG_i) = |R(\cG_i)| - |B(\cG_i)| = |R(\cG_i)| - \bigl\lceil |R(\cG_i)|/3 \bigr\rceil, 
			\]
			and hence from \Cref{lem:discrepancy} and the fact $|R(\cG_i)| \ge 2$ we deduce that 
			\begin{equation} \label{eq:gs_rrbbound}
				\Nrrb(\cG_i) \ge |B(\cG_i)| \cdot \disc(\cG_i) \ge \frac{2}{9} |R(\cG_i)|^2. 
			\end{equation}
			Here we applied the fact that $\lc x/3 \rc \cdot \bigl( x - \lc x/3 \rc \bigr) \ge 2x^2/9$ holds for any integer $x \ge 2$. 
			\vspace{-0.5em}
			\item Observe that the sector $\cT$ contains all red points in $R(\cH)$ that are not in $\cG_1 \cup \dots \cup \cG_k$, and $|B(\cT)| < |R(\cT)|/3$. Reflect $R \cup B$ over the $y$-axis and conduct the same algorithm to generate good sectors $\cG_1', \dots, \cG_{k'}'$. Define $\cT'$ in the same way (with respect to the new family of good sectors). Then we reflect everything over the $y$-axis again to work on the original picture. See \hyperlink{figthree}{Figure~3} for an illustration of the same example as in \hyperlink{figtwo}{Figure~2}. 
		\end{itemize}
		
		\vspace{0.5em}
		\begin{center}
			\begin{tikzpicture}[line cap=round,line join=round,>=triangle 45,x=2cm,y=2cm,scale = 0.9]
				\clip(-2,-0.7) rectangle (2,2);
				\fill[line width=0pt,color=qqqqff,fill=qqqqff,fill opacity=0.1] (0,0) -- (4.617395575649076,0) -- (4.492019706117622,1.2227409713247717) -- cycle;
				\fill[line width=0pt,color=ccqqqq,fill=ccqqqq,fill opacity=0.1] (0,0) -- (4.492019706117622,1.2227409713247717) -- (4.405823795814748,2.398139748182151) -- cycle;
				\fill[line width=0pt,color=ccqqqq,fill=ccqqqq,fill opacity=0.1] (0,0) -- (2.1725661197857247,2.954495169227977) -- (0.02166145868131356,3.2899724801546206) -- cycle;
				\fill[line width=0pt,color=qqqqff,fill=qqqqff,fill opacity=0.1] (0,0) -- (-0.8611262409352076,3.2290905698362398) -- (-1.3938429562210393,3.046444838881098) -- cycle;
				\fill[line width=0pt,color=qqqqff,fill=qqqqff,fill opacity=0.1] (0,0) -- (-3.735771731884041,1.324608865319078) -- (-1.3938429562210393,3.046444838881098) -- cycle;
				\fill[line width=0pt,color=qqqqff,fill=qqqqff,fill opacity=0.1] (0,0) -- (-0.8611262409352076,3.2290905698362398) -- (0.02166145868131356,3.2899724801546206) -- cycle;
				\fill[line width=0pt,color=qqqqff,fill=qqqqff,fill opacity=0.1] (0,0) -- (2.1725661197857247,2.954495169227977) -- (4.405823795814748,2.398139748182151) -- cycle;
				\fill[line width=0pt,color=ccqqqq,fill=ccqqqq,fill opacity=0.1] (0,0) -- (-3.4771840009754174,0) -- (-3.735771731884041,1.324608865319078) -- cycle;
				
				\draw [line width=0.8pt,domain=-4.026985577491008:4.46604091192311,dashed] plot(\x,{(-0-0*\x)/0.5});
				\draw [line width=0.8pt,domain=0:4.46604091192311] plot(\x,{(-0--2.398139748182151*\x)/4.405823795814748});
				\draw [line width=0.8pt,domain=0:1] plot(\x,{(-0--3.2899724801546206*\x)/0.02166145868131356});
				
				\draw [line width=0.8pt,domain=-4.026985577491008:0] plot(\x,{(-0--3.046444838881098*\x)/-1.3938429562210393});
				\draw [line width=0.8pt,domain=-4.026985577491008:0] plot(\x,{(-0--2.7800864812381825*\x)/-3.4257767130969974});
				\draw [line width=0.8pt,domain=-4.026985577491008:0] plot(\x,{(-0--3.8759608669690344*\x)/-0.4730040626555302});
				\draw [line width=0.8pt,domain=0:4.46604091192311] plot(\x,{(-0--2.9246810182443364*\x)/4.1388006439618135});
				
				\draw [line width=0.8pt,domain=0:0.2] plot({-1.0947051280427746+\x*(cos(110))},{0+\x*(sin(110))});
				\draw [line width=0.8pt,domain=0:0.2] plot({-0.8500234991900003+\x*(cos(205))},{0.6898111105135488+\x*(sin(205))});
				\draw [line width=0.8pt,domain=0:0.2] plot({0.9614985018924833+\x*(cos(340))},{0.5233545148574316+\x*(sin(340))});
				\draw [line width=0.8pt,domain=0:0.2] plot({0.8940157100276518+\x*(cos(100))},{0.6317556707991636+\x*(sin(100))});
				\draw [line width=0.8pt,domain=0:0.2] plot({-0.4554539113344601+\x*(cos(90))},{0.9954602212109177+\x*(sin(90))});
				\draw [line width=0.8pt,domain=0:0.2] plot({0.007207475466635341+\x*(cos(20))},{1.0946814009841153+\x*(sin(20))});
				\draw [line width=0.8pt,domain=0:0.2] plot({1.0947051280427746+\x*(cos(60))},{0+\x*(sin(60))});
				\draw [line width=0.8pt,domain=0:0.2] plot({-0.13260887587739487+\x*(cos(120))},{1.0866435493765567+\x*(sin(120))});
				
				\draw [shift={(0,0)},line width=0.8pt]  plot[domain=0:0.49846521055941184,variable=\t]({1*1.0947051280427746*cos(\t r)+0*1.0947051280427746*sin(\t r)},{0*1.0947051280427746*cos(\t r)+1*1.0947051280427746*sin(\t r)});
				\draw [shift={(0,0)},line width=0.8pt]  plot[domain=0.6151747263955964:1.5642123368646565,variable=\t]({1*1.0947051280427746*cos(\t r)+0*1.0947051280427746*sin(\t r)},{0*1.0947051280427746*cos(\t r)+1*1.0947051280427746*sin(\t r)});
				\draw [shift={(0,0)},line width=0.8pt]  plot[domain=1.6922311788190205:1.999895357801293,variable=\t]({1*1.0947051280427746*cos(\t r)+0*1.0947051280427746*sin(\t r)},{0*1.0947051280427746*cos(\t r)+1*1.0947051280427746*sin(\t r)});
				\draw [shift={(0,0)},line width=0.8pt]  plot[domain=2.459866659096932:3.141592653589793,variable=\t]({1*1.0947051280427749*cos(\t r)+0*1.0947051280427749*sin(\t r)},{0*1.0947051280427749*cos(\t r)+1*1.0947051280427749*sin(\t r)});
				
				\draw [fill=ffqqtt, draw=none] (0,0) circle (2pt);
				\draw[color=black] (0,-0.15) node {\footnotesize $r_0$};
				\draw[color=black] (1.2,0.32) node {\footnotesize $\cT'$};
				\draw[color=black] (0.61,1.09) node {\footnotesize $\cG'_3$};
				\draw[color=black] (-0.325,1.2) node {\footnotesize $\cG'_2$};
				\draw[color=black] (-1.16,0.42) node {\footnotesize $\cG'_1$};
				
				\draw [fill=ffqqtt, draw=none] (-1.5917091647557768,0) circle (2pt);
				\draw [fill=ffqqtt, draw=none] (-0.778822214656253,0.2761503871439254) circle (2pt);
				\draw [fill=ffqqtt, draw=none] (-0.7501682436628383,1.6396009062564258) circle (2pt);
				\draw [fill=ffqqtt, draw=none] (0.004752251305354271,0.7217785396364786) circle (2pt);
				\draw [fill=ffqqtt, draw=none] (1.0619597233890699,1.4441700273670546) circle (2pt);
				\draw [fill=ffqqtt, draw=none] (1.599970966599487,0.8708823023255056) circle (2pt);
				\draw [fill=ffqqtt, draw=none] (1.7676837407854369,0.48116871153957375) circle (2pt);
				\draw [fill=ffqqtt, draw=none] (-0.8687364797250052,0.12030802420392708) circle (2pt);
				\draw [fill=ffqqtt, draw=none] (-1.8200163284497048,0.32578447152846185) circle (2pt);
				\draw [fill=qqqqff, draw=none] (-1.4740963834589054,1.2950521573926819) circle (2pt);
				\draw [fill=qqqqff, draw=none] (-1.667811552653753,0.8661114256040903) circle (2pt);
				\draw [fill=qqqqff, draw=none] (-1.3938429562210393,1.4216588572593138) circle (2pt);
				\draw [fill=qqqqff, draw=none] (-0.6632600324004702,0.6758554558591506) circle (2pt);
				\draw [fill=qqqqff, draw=none] (-0.9905003003617667,0.9269933359224709) circle (2pt);
				\draw [fill=ffqqtt, draw=none] (0.32953020972312486,1.0058630833803737) circle (2pt);
				\draw [fill=ffqqtt, draw=none] (0.6034988061558383,1.622292425353978) circle (2pt);
				\draw [fill=ffqqtt, draw=none] (0.7902955764508708,0.341004949108057) circle (2pt);
				\draw [fill=qqqqff, draw=none] (1.5547786548805378,0.1915675328720319) circle (2pt);
				\draw [fill=qqqqff, draw=none] (-0.22921946580604066,1.878304542298994) circle (2pt);
				\draw [fill=qqqqff, draw=none] (1.8171696698845659,1.2841018685435752) circle (2pt);
				\draw [fill=qqqqff, draw=none] (1.5020203663260405,1.3780152913897212) circle (2pt);
				\draw [fill=qqqqff, draw=none] (-0.9770850399750475,1.169012138532197) circle (2pt);
				\draw [fill=qqqqff, draw=none] (-1.3758796755478557,1.1165539397758957) circle (2pt);
				\draw [fill=qqqqff, draw=none] (-0.5,0.5) circle (2pt);
				\draw [fill=qqqqff, draw=none] (-1.1772970706567991,1.2798316798130869) circle (2pt);
				\draw [fill=qqqqff, draw=none] (-1.6186909204650597,1.4624774107682288) circle (2pt);
				\draw [fill=qqqqff, draw=none] (-1.519757816197691,1.7212255296213468) circle (2pt);
				\draw [fill=qqqqff, draw=none] (-1.801336651420202,1.6679538580927638) circle (2pt);
				\draw [fill=ffqqtt, draw=none] (1.3873534015049906,0.5644692335721136) circle (2pt);
				
				\node at (0, -0.55) {\textbf{\hypertarget{figthree}{Figure 3}:} Generating $\cG'_1, \dots, \cG'_{k'}$ and $\cT'$. };
				
			\end{tikzpicture}
		\end{center}
		
		\begin{itemize}[resume]
			\item[] Similarly, $|B(\cT')| < |R(\cT')|/3$. We claim that $\cT \cap \cT' = \{r_0\}$. If not, then $\cT \cup \cT' = \cH$, and so 
			\[
			|B(\cH)| \le |B(\cT)| + |B(\cT')| < \frac{|R(\cT)|}{3} + \frac{|R(\cT')|}{3} \le \frac{2|R(\cH)|}{3} \le |R(\cH)| - 1, 
			\]
			where in the last inequality we applied $|R(\cH)| \ge 3$. This contradicts $|B(\cH)| \ge |R(\cH)| - 1$. It then follows that $|R(\cT)| + |R(\cT')| \le |R(\cH)| + 1$, as $r_0$ is counted twice. Assume without loss of generality that $|R(\cT)| \le \bigl( |R(\cH)| + 1 \bigr) / 2$. Since $|R(\cH)| \ge 2$, we have $|R(\cT)| < |R(\cH)|$, and hence $k \ge 1$. By counting $r_0$ to the multiplicity $k+1$, we derive that
			\[
			|R(\cG_1)| + \dots + |R(\cG_k)| + |R(\cT)| = |R(\cH)| + k. 
			\]
			This together with $R(\cT) \le \bigl( |R(\cH)| + 1 \bigr) / 2, \, k \ge 1$, and $|R(\cH)| = \bigl\lceil (n+1)/2 \bigr\rceil$ implies that
			\begin{equation} \label{eq:gs_large}
				\sum_{i=1}^k |R(\cG_i)| \ge \frac{|R(\cH)|-1}{2} + k \ge \frac{n}{4} \implies \frac{1}{k}\sum_{i=1}^k |R(\cG_i)| \ge \frac{n}{4k}. 
			\end{equation}
		\end{itemize}
		By putting all the ingredients from above together, we conclude that
		\begin{align*}
			\Nrrb(\R^2) \ge \sum_{i=1}^k \Nrrb(\cG_i) &\ge \sum_{i=1}^k \frac{2}{9} |R(\cG_i)|^2 &&\text{by \eqref{eq:gs_rrbbound}} \\
			&\ge k \cdot \frac{2}{9} \Bigl( \frac{n}{4k} \Bigr)^2 &&\text{by convexity and \eqref{eq:gs_large}} \\
			&\ge \frac{n^2}{72\sqrt{n}} &&\text{by \eqref{eq:gs_stepnumber}} \\
			&= \Omega(n^{3/2}). 
		\end{align*}
		The proof of \Cref{thm:RRB_3/2} is complete. 
	\end{proof}
	
	\section{Concluding remarks}
	
	We do not have any construction of $n$ red and $n$ blue points in the plane containing subquadratically many empty red-red-blue triangles. In the literature, Horton sets are often good candidates for such extremal constructions. However, it is not hard to see that every balanced bi-coloring of an $n$-point Horton set contains $\Omega(n^2)$ empty red-red-blue triangles. We state the following conjecture: 
	
	\begin{conjecture}
		There exist $\Omega(n^2)$ empty red-red-blue triangles in $R \cup B$. 
	\end{conjecture}
	
	\section*{Acknowledgments}
	
	The second author benefited from discussions with J\'{a}nos Pach at an early stage of this work. We are grateful to J\'{a}nos Pach for suggesting this problem and for useful comments. 
	
	After posting an earlier arXiv version of this paper, we learned that our main result was proved independently by D\'{\i}az-ba\~{n}ez, Fabila-Monroy, and Urrutia \cite{diaz-banez_fabila-monroy_urrutia}. We are grateful to Oriol Sol\'{e} Pi for pointing \cite{diaz-banez_fabila-monroy_urrutia} to us. We also thank Adrian Dumitrescu for telling us about \cite{aichholzer_etal_5hole}. 
	
	\bibliographystyle{plain}
	\bibliography{RRB}
	
\end{document}